\pdfoutput=1
\documentclass[11pt,a4paper]{article}

\usepackage[english]{babel}

\usepackage[numbers,square,sort&compress]{natbib}

\usepackage[
    a4paper,
    top=2.5cm,
    bottom=2.5cm,
    left=3cm,
    right=3cm
]{geometry}

\usepackage{setspace}
\usepackage{amsmath}
\usepackage{amssymb}
\usepackage{amsthm}
\usepackage{mathtools}

\usepackage{graphicx}
\usepackage{caption}

\usepackage{algorithm}
\usepackage{algpseudocode}

\usepackage{booktabs}
\usepackage{makecell}
\usepackage{array}

\usepackage[colorlinks=true, allcolors=blue]{hyperref}

\theoremstyle{definition}

\newtheorem{definition}{Definition}[section]

\theoremstyle{plain}

\newtheorem{proposition}{Proposition}[section]

\title{Mean Field Games in Dynamic Ticket Pricing with Overselling and Cancellations}

\author{
Yuzhou Chen\thanks{Faculty of Electrical Engineering, Mathematics and
Computer Science, Delft University of Technology.
Email: \texttt{yuzhouchen028@outlook.com}.}
\and
Yumiharu Nakano\thanks{Department of Mathematical and Computing Science,
Institute of Science Tokyo; Future Intelligence Visionary Initiative,
Institute of Science Tokyo. Corresponding author.
Email: \texttt{nakano@comp.isct.ac.jp}.}
}
\date{}

\begin{document}

\maketitle

\begin{abstract}
We extend a finite-state Mean Field Game model of dynamic ticket pricing to
allow overselling, stochastic cancellations, and terminal penalties for
unresolved negative inventory. Cancellations generate upward inventory
transitions, producing a bidirectional inventory process coupled with firms'
pricing decisions. We formulate the representative-firm control problem and
the population consistency conditions, and prove existence of a selected
Markovian mean field equilibrium under the stated assumptions, including a
compact price set and an explicit tie-breaking rule. Numerical experiments show
that the quadratic terminal penalty can generate an endogenous soft-boundary
mechanism that discourages progressively deeper overselling. Under a
low-cancellation stress case, linear penalties more strongly reduce the
prevalence of terminal overselling, whereas quadratic penalties reduce its
conditional depth and lower-tail severity. A time-varying cancellation
experiment suggests that anticipated future inventory recovery can support
lower early prices and greater temporary overbooking exposure. Under the
adopted demand normalization, stronger mean-field coupling is associated with
lower equilibrium prices and slower inventory depletion. Together, these
results provide a tractable framework for analyzing pricing and overbooking
incentives in competitive ticket markets with cancellations.
\end{abstract}


\section{Introduction}
\label{sec:introduction}

Airline revenue management concerns the dynamic pricing and allocation of
perishable capacity under uncertain demand
\citep{talluri2004revenue,phillips2005pricing,klein2020review}.
A central complication is overbooking: airlines may accept bookings beyond
physical capacity to hedge against cancellations and no-shows, but unresolved
excess bookings at departure may generate compensation, rebooking, and other
operational costs
\citep{rothstein1971airline,subramanian1999yield,aydin2013single}.

Modeling these decisions in a competitive market becomes increasingly difficult
as the number of firms grows. Mean Field Game (MFG) theory provides a tractable
large population framework in which a representative firm interacts with
the aggregate market distribution. In dynamic ticket pricing,
\citet{aydin2026fare} develop a finite-state MFG in which firms control
sales intensity through prices and interact through the mean field. Their
baseline model, however, imposes a strict absorbing boundary at zero inventory,
so that overselling and cancellation-driven inventory recovery are excluded.

\citet{aydin2026fare} identify overselling with penalty functions
and cancellation dynamics as natural directions for extending the
strict-capacity framework. We develop these extensions by allowing negative
inventory states, introducing stochastic cancellation-driven upward transitions,
and imposing terminal penalties for unresolved overbooking. The resulting model
makes four contributions. First, it formulates a bidirectional finite-state MFG
for overselling and cancellation-driven inventory recovery and defines
equilibrium through the population state distribution and mean quoted price.
Second, under the stated assumptions, including a compact price set and an explicit
tie-breaking rule, it proves existence of a selected Markovian mean field
equilibrium while separating existence from global uniqueness. Third, it shows
numerically that the quadratic terminal penalty can generate an endogenous
soft-boundary mechanism that discourages deeper overselling. Fourth, it studies
cancellation timing and mean-field coupling: the time-varying experiment
suggests that anticipated future recovery can support lower early prices and
greater temporary overbooking exposure, while, under the adopted demand
normalization, stronger coupling is associated with lower prices and slower
inventory depletion.

The selected equilibrium is characterized by a coupled system consisting of
the backward Hamilton--Jacobi--Bellman (HJB) equations and the forward
Kolmogorov equations, and is computed numerically using a fixed-point
iteration scheme. The equilibrium concept is the non-atomic MFG
equilibrium; no finite-player Nash approximation result is claimed.
The baseline formulation serves as the numerical starting point for the extended
model and is summarized in Appendix~\ref{app:baseline}.

\section{Related Literature}
\label{sec:literature}

Mean Field Games (MFGs) provide a tractable framework for large-population
dynamic games by replacing the direct interaction among many agents with the
interaction between a representative agent and an aggregate population
distribution \citep{lasry2007mean,huang2006large,carmona2018probabilistic}.
While much of the early literature considers continuous state spaces,
finite-state MFGs are particularly relevant for discrete inventory problems
\citep{gomes2013finite,cecchin2020probabilistic,gomes2021numerical}. Extended
finite-state MFGs may also involve interaction through the distribution of
controls rather than only through the state distribution
\citep{carmona2021extended,garcia2025control}. Computational and learning-based
approaches to MFGs have also developed rapidly
\citep{lauriere2022learning,magnino2025solving}. In dynamic pricing,
\citet{aydin2026fare} formulate a finite-state MFG in which firms control
sales intensity through prices while interacting through the aggregate market.

The present work is also related to the broader revenue-management literature
on dynamic pricing and stochastic demand. Classical models study the allocation
and pricing of perishable capacity
\citep{talluri2004revenue,phillips2005pricing,klein2020review}, including
expected marginal seat revenue methods \citep{belobaba1989application} and
stochastic intensity control \citep{gallego1994optimal}. Review articles discuss
multiple-product and competitive extensions
\citep{bitran2003overview,chen2015recent}. \citet{gallego2014dynamic} directly
study dynamic pricing under competition, while mean-field approaches provide an
alternative way to represent strategic interaction in large markets
\citep{paraschiv2021dynamic,aydin2026fare}.

A separate revenue-management literature studies overbooking, cancellations,
and no-shows. Overbooking allows firms to accept reservations beyond physical
capacity in order to hedge against uncertain booking realization
\citep{rothstein1971airline,talluri2004revenue,subramanian1999yield}.
Dynamic formulations incorporate cancellation risk and time-to-departure effects
\citep{subramanian1999yield,aydin2013single,alavifard2019optimal,
dai2019network}. Customer choice with cancellations and overbooking is treated
by \citet{sierag2015revenue}. Cancellations naturally generate
bidirectional inventory dynamics, since sales decrease available inventory
whereas cancellations restore it. By contrast, the finite-state ticket-pricing
MFG of \citet{aydin2026fare} uses a strict zero-inventory boundary.
The present paper combines these strands by introducing negative inventory
states, cancellation-driven upward transitions, and terminal overbooking
penalties within the finite-state MFG framework.

\section{Mean Field Game Framework with Overselling and Cancellations}
\label{sec:extension}
% ============================================================

This section develops the Mean Field Game framework used in the remainder of
the paper. We begin with a concise summary of the finite-state ticket-pricing
model of \citet{aydin2026fare}, since its pricing mechanism,
demand specification, and mean-field interaction are retained in our analysis.
We then extend the inventory dynamics to allow negative inventory states and
cancellation-driven recovery, and introduce terminal penalties for unresolved
overbooking.

The purpose of retaining the baseline structure is to isolate the effects of
these extensions. In particular, the enlarged state space changes the inventory
process from one-sided depletion to bidirectional dynamics, while the terminal
penalty introduces an additional intertemporal trade-off between current sales,
future cancellation recovery, and overbooking risk at departure.

% ============================================================
\subsection{Baseline framework}
\label{sec:baseline}
% ============================================================

The baseline model describes a continuum of firms that dynamically control
their sales intensity through prices. Inventory takes values in
\[
\{0,1,\dots,K\},
\]
and decreases by one unit whenever a sale occurs. Once inventory reaches
$k=0$, further sales are prohibited, making zero inventory an absorbing sales
boundary.

For a given mean field, the representative firm's optimization problem is
characterized by a backward Hamilton--Jacobi--Bellman (HJB) equation, while
the population distribution evolves according to a forward Kolmogorov
equation. Market consistency between individual pricing decisions and the
aggregate distribution determines the mean field equilibrium.

The baseline formulation provides the numerical starting point for the
extended model. Its core mathematical structure and numerical implementation
are summarized in Appendix~\ref{app:baseline}.

% ============================================================
\subsection{Overselling and cancellations}
\label{sec:overselling-cancellations}
% ============================================================

If sales are restricted to the physical capacity $K$, cancellations may leave
seats unused at departure. Overbooking provides a natural hedge against this
risk by allowing firms to accept bookings beyond physical capacity, but
unresolved excess bookings may generate compensation, rebooking, and other
operational costs
\citep{talluri2004revenue,subramanian1999yield,aydin2013single}.

To represent this trade-off, we allow the inventory state to become negative.
Negative inventory represents bookings accepted beyond physical capacity,
while cancellations generate upward transitions that may restore inventory
before departure.

\begin{definition}[Extended State Space]
Let $K$ denote physical capacity and let $M>0$ denote the maximum admissible
overselling depth. The extended state space is
\begin{equation}
    \mathcal{X}_{\mathrm{ext}}
    =
    \{-M,-M+1,\dots,0,1,\dots,K\}.
    \label{eq:extended_state_space}
\end{equation}
\end{definition}

A state $k<0$ indicates that the firm has accepted $|k|$ bookings beyond
physical capacity. The state $k=0$ is therefore no longer an absorbing sales
boundary. Sales may continue until the lower boundary $k=-M$ is reached,
whereas cancellations move inventory upward from $k$ to $k+1$.

The parameter $M$ specifies the lower boundary of the finite inventory state
space. It is treated as a stylized modeling parameter rather than as an
empirically calibrated airline booking limit. Its benchmark value is reported
in Section~\ref{sec:results}.

The firm therefore balances the immediate revenue from additional sales
against two opposing forces: future cancellations may restore inventory,
whereas unresolved negative inventory at departure generates a terminal
penalty.

% ============================================================
\subsection{Extended HJB equation}
\label{sec:extended-hjb}
% ============================================================

Allowing negative inventory and cancellations modifies the dynamic programming
problem in two ways. Sales may continue for $k<0$ subject to the lower boundary
$-M$, while cancellations generate upward transitions from $k$ to $k+1$.

% ------------------------------------------------------------
\subsubsection{Demand specification and competitive coupling}
\label{sec:demand-competition}
% ------------------------------------------------------------

We retain the linear sales-intensity specification of the baseline model,
\begin{equation}
    \lambda(p,\bar{p},m)
    =
    \left(
        a(m)-p+c(m)\bar{p}
    \right)^+,
    \label{eq:extended_linear_intensity}
\end{equation}
where
\begin{equation}
    a(m)
    =
    \frac{1}{1+\epsilon\eta_{\mathrm{ext}}(t)},
    \qquad
    c(m)
    =
    \frac{\epsilon\eta_{\mathrm{ext}}(t)}
    {1+\epsilon\eta_{\mathrm{ext}}(t)}.
    \label{eq:competition_coefficients}
\end{equation}

This specification is inherited from the baseline model so that overselling
and cancellations can be studied without changing the underlying demand
structure. The component of sales intensity independent of the firm's own
price is
\begin{equation}
    b(t)
    :=
    a(m)+c(m)\bar{p}(t)
    =
    \frac{
        1+\epsilon\eta_{\mathrm{ext}}(t)\bar{p}(t)
    }{
        1+\epsilon\eta_{\mathrm{ext}}(t)
    }.
    \label{eq:effective_demand_component}
\end{equation}

Hence, changing $\epsilon$ affects both the weight placed on the mean quoted
price and the effective demand level. In particular, when $\bar{p}(t)<1$,
increasing $\epsilon$ lowers $b(t)$, holding the endogenous market quantities
fixed. This normalization is therefore important when interpreting the
comparison in Section~\ref{sec:competition-impact}; the effect of $\epsilon$
should be viewed as a model-specific comparative static rather than as a
general prediction about market concentration and overbooking.

For the existence analysis, admissible prices take values in the compact set
\begin{equation}
    \mathcal A=[0,P_{\max}],
    \qquad P_{\max}\geq1.
    \label{eq:compact_action_set}
\end{equation}
The cap is a technical compactness assumption. It is chosen above the selected
prices attained in the reported computations and is therefore nonbinding in
the numerical results. Indeed, if $0\leq\bar p\leq P_{\max}$, then
$0\leq b(t)\leq P_{\max}$.

% ------------------------------------------------------------
\subsubsection{Cancellation dynamics}
\label{sec:cancellation-dynamics}
% ------------------------------------------------------------

During the selling horizon $t\in[0,T)$, a sale decreases inventory by one unit,
whereas a cancellation restores one unit. Following revenue management models
that treat cancellations as stochastic components of the booking process
\citep{subramanian1999yield,sierag2015revenue,dai2019network},
we model cancellations as exogenous Poisson jumps.

If $X_t=k$, then $K-k$ bookings are outstanding relative to the initial
capacity. We assume that each outstanding ticket is canceled independently at
rate $\gamma>0$. The resulting aggregate cancellation intensity is
\begin{equation}
    \mu(k)
    =
    \gamma(K-k)^+.
    \label{eq:constant_cancellation_intensity}
\end{equation}

Thus, the cancellation intensity increases with the number of outstanding
bookings. Cancellation behavior is treated as exogenous and does not depend on
the firm's current price or on individual customer characteristics.

When a cancellation occurs, inventory moves from $k$ to $k+1$, and the firm
incurs a cancellation-related cost $R\geq0$. The parameters $\gamma$ and $R$
are treated as normalized model inputs rather than empirically calibrated
quantities. Their benchmark values are reported in Section~\ref{sec:results}.

To examine the effect of anticipated future inventory recovery, we also
consider a time-dependent cancellation rate,
\begin{equation}
    \mu(k,t)
    =
    \gamma(t)(K-k)^+,
    \label{eq:time_dependent_cancellation_intensity}
\end{equation}
where
\begin{equation}
    \gamma(t)
    =
    \gamma_{\mathrm{base}}
    +
    \gamma_{\mathrm{surge}}
    \left(
        \frac{t}{T}
    \right)^\beta,
    \qquad
    \beta>1.
    \label{eq:gamma_time_dependent}
\end{equation}

Here, $\gamma_{\mathrm{base}}$ is the baseline per-ticket cancellation rate,
$\gamma_{\mathrm{surge}}$ controls its increase toward departure, and $\beta$
determines how strongly this increase is concentrated near the terminal date.
This specification is used for scenario analysis rather than as an empirical
description of airline cancellation behavior.

Let $N^s$ and $N^c$ be simple counting processes with predictable intensities
\[
 \mathbf 1_{\{X_{t-}>-M\}}
 \lambda(p_t,\bar p(t),m(t))
 \quad\text{and}\quad
 \mathbf 1_{\{X_{t-}<K\}}\mu(X_{t-},t),
\]
respectively. We assume that simultaneous sales and cancellation jumps have
probability $o(h)$ over an interval of length $h$. The extended inventory
process satisfies
\begin{equation}
    dX_t=-dN_t^s+dN_t^c.
    \label{eq:controlled_inventory_dynamics}
\end{equation}
Since $\mu(K,t)=0$, the upper-boundary indicator in the cancellation intensity
may equivalently be omitted. This point-process construction follows the
standard intensity-based framework; see, for example,
\citet{bremaud1981point}. The associated finite-state control problem is a
continuous-time Markov decision problem of the type treated in
\citet{guo2009ctmdp}.

For a prescribed deterministic population flow
$(m(t),\bar p(t))_{t\in[0,T]}$, an admissible Markov feedback is a measurable map
\[
    \alpha:[0,T]\times\mathcal X_{\mathrm{ext}}\longrightarrow\mathcal A.
\]
We set $\alpha(t,-M)=0$, since the quoted price at the lower boundary affects
neither transitions nor the mean quoted price. The controlled generator is
\begin{equation}
\begin{aligned}
(\mathcal L_t^{\alpha;m,\bar p}\varphi)(k)
={}&\mathbf 1_{\{k>-M\}}
 \lambda(\alpha(t,k),\bar p(t),m(t))
 [\varphi(k-1)-\varphi(k)]\\
&+\mathbf 1_{\{k<K\}}\mu(k,t)
 [\varphi(k+1)-\varphi(k)].
\end{aligned}
\label{eq:controlled_generator}
\end{equation}
Let $\Psi:\mathcal X_{\mathrm{ext}}\to[0,\infty)$ be the terminal penalty
function, with $\Psi(k)=0$ for $k\geq0$; its specifications for $k<0$ are given
in Section~\ref{sec:penalty-design}. Given $(m,\bar p)$, the representative
firm maximizes
\begin{equation}
J_{t,k}(\alpha;m,\bar p)
=\mathbb E_{t,k}\!\left[
 \int_t^T e^{-r(s-t)}
 \bigl(\alpha(s,X_{s-})\,dN_s^s-R\,dN_s^c\bigr)
 -e^{-r(T-t)}\Psi(X_T)
\right],
\label{eq:representative_objective}
\end{equation}
and its value is
\begin{equation}
    V_k(t;m,\bar p)=\sup_{\alpha}J_{t,k}(\alpha;m,\bar p).
    \label{eq:representative_value}
\end{equation}
Here $R$ is a normalized cost per cancellation; it is not modeled as a refund
of the historical transaction price.

For $0<h\leq T-t$, the dynamic programming principle (DPP) over
$[t,t+h]$ is
\begin{equation}
\begin{aligned}
V_k(t;m,\bar p)=\sup_{\alpha}\mathbb E_{t,k}\Bigg[&
 \int_t^{t+h}e^{-r(s-t)}\alpha(s,X_{s-})\,dN_s^s
 -\int_t^{t+h}e^{-r(s-t)}R\,dN_s^c\\
&+e^{-rh}V_{X_{t+h}}(t+h;m,\bar p)\Bigg].
\end{aligned}
\label{eq:representative_dpp}
\end{equation}
Here the supremum is taken over admissible controls on $[t,t+h]$, after
which the continuation value is used. For an interior state, set
$p=\alpha(t,k)$. The sale, cancellation, and no-jump probabilities are
$\lambda(p,\bar p,m)h+o(h)$, $\mu(k,t)h+o(h)$, and
$1-(\lambda(p,\bar p,m)+\mu(k,t))h+o(h)$, respectively. We substitute these
expressions into Equation~\eqref{eq:representative_dpp}, use
$e^{-rh}V_k(t+h)=V_k(t)+h(\frac{dV_k}{dt}-rV_k)+o(h)$ in the no-jump term
and zeroth-order continuation values in the jump terms, subtract $V_k(t)$,
divide by $h$, and let $h\downarrow0$. This yields, for
$k\in\{-M+1,\dots,K-1\}$, the HJB equation
\begin{equation}
    \frac{dV_k}{dt}
    -
    rV_k
    +
    \sup_{p\in\mathcal A}
    \left\{
        \lambda(p,\bar{p},m)
        \bigl(
            p-\Delta V_k(t)
        \bigr)
    \right\}
    +
    \mu(k,t)
    \left[
        \Delta V_{k+1}(t)-R
    \right]
    =
    0,
    \label{eq:extended_hjb}
\end{equation}
where
\begin{equation}
    \Delta V_k(t)
    =
    V_k(t)-V_{k-1}(t).
\end{equation}
In particular, the value increment associated with a cancellation from $k$ to
$k+1$ is $V_{k+1}(t)-V_k(t)=\Delta V_{k+1}(t)$.

The sales term represents the immediate revenue from an additional sale net of
the marginal value of inventory. The cancellation term represents the expected
change in continuation value associated with an upward transition, net of the
cancellation-related cost $R$. Future cancellations may therefore increase the
continuation value of negative inventory states and affect current pricing
decisions.

To describe the optimizer, write $C_k(t)=\Delta V_k(t)$ and retain the notation
$b(t)$ from Equation~\eqref{eq:effective_demand_component}. Since
$0<b(t)\leq P_{\max}$, the maximizers of the sales Hamiltonian over the compact
admissible set are
\begin{equation}
\operatorname*{arg\,max}_{p\in\mathcal A}
 \bigl(b(t)-p\bigr)^+(p-C_k)
=
\begin{cases}
 \{0\}, & C_k<-b(t),\\[1mm]
 \{(b(t)+C_k)/2\}, & -b(t)\leq C_k<b(t),\\[1mm]
 [b(t),P_{\max}], & C_k\geq b(t).
\end{cases}
\label{eq:price_argmax_set}
\end{equation}
When $C_k(t)\geq b(t)$, every admissible price in $[b(t),P_{\max}]$ is optimal.
This optimizer is non-unique whenever $b(t)<P_{\max}$.
Throughout the paper we adopt the continuous tie-breaking selection
\begin{equation}
 p_k^*(t)=
 \min\left\{
 P_{\max},
 \frac12\bigl(b(t)+C_k(t)\bigr)^+
 \right\},
 \qquad k>-M,
 \label{eq:selected_price}
\end{equation}
When the cap is nonbinding, the minimum does not alter the pricing formula used
in the numerical implementation. Consequently, all reported state-dependent
and mean quoted prices refer to this selected feedback rule.

Although the quoted price can be non-unique, the induced optimal sales
intensity is uniquely determined:
\begin{equation}
 \lambda_k^*(t)=
 \begin{cases}
 b(t), & C_k(t)<-b(t),\\[1mm]
 \dfrac{b(t)-C_k(t)}{2}, & -b(t)\leq C_k(t)<b(t),\\[2mm]
 0, & C_k(t)\geq b(t).
 \end{cases}
 \label{eq:unique_sales_intensity}
\end{equation}

% ------------------------------------------------------------
\subsubsection{Terminal penalties}
\label{sec:terminal-penalties}
% ------------------------------------------------------------

Unresolved negative inventory at the terminal date generates an overbooking
penalty, while nonnegative terminal inventory has zero salvage value. Let
$\Psi(k)\geq0$ denote the penalty associated with ending at state $k<0$. The
terminal condition is
\begin{equation}
    V_k(T)
    =
    \begin{cases}
        0,
        & k\geq0,
        \\[0.3em]
        -\Psi(k),
        & k<0.
    \end{cases}
    \label{eq:terminal_penalty_condition}
\end{equation}

The function $\Psi$ provides a reduced-form representation of compensation,
rebooking, operational disruption, and other costs associated with unresolved
overbooking. Because the HJB equation is solved backward in time, the terminal
penalty also influences pricing before departure through the continuation
value of negative inventory states.

% ------------------------------------------------------------
\subsubsection{Boundary conditions}
\label{sec:extended-boundaries}
% ------------------------------------------------------------

At the lower boundary $k=-M$, further sales are not permitted, although
cancellations may still move the firm upward to $-M+1$. The corresponding HJB
equation is
\begin{equation}
    \frac{dV_{-M}}{dt}
    -
    rV_{-M}
    +
    \mu(-M,t)
    \left[
        \Delta V_{-M+1}(t)-R
    \right]
    =
    0.
    \label{eq:hjb_lower_boundary}
\end{equation}

Thus, $k=-M$ is a lower sales boundary rather than a fully absorbing state. At
the upper boundary $k=K$, no tickets have been sold relative to physical
capacity and
\[
\mu(K,t)=0,
\]
so the HJB equation contains the sales term but no upward cancellation term:
\begin{equation}
    \frac{dV_K}{dt}-rV_K
    +\sup_{p\in\mathcal A}
    \left\{
        \lambda(p,\bar p,m)\bigl(p-\Delta V_K(t)\bigr)
    \right\}=0.
    \label{eq:hjb_upper_boundary}
\end{equation}

% ============================================================
\subsection{Penalty functions}
\label{sec:penalty-design}
% ============================================================

Terminal penalties represent the cost of unresolved overbooking at departure
\citep{rothstein1971airline,talluri2004revenue,aydin2013single}. We compare
linear and quadratic specifications as reduced-form representations of constant
and increasing marginal penalty profiles, respectively.

% ------------------------------------------------------------
\subsubsection{Linear penalty}
% ------------------------------------------------------------

For $c_1>0$, define
\begin{equation}
    \Psi_{\mathrm{lin}}(k)
    =
    c_1|k|,
    \qquad
    k<0.
    \label{eq:linear_penalty}
\end{equation}

For $k<0$, the terminal condition gives
$V_k(T)=-c_1|k|$, and hence
\begin{align}
    \Delta V_k(T)
    &=
    V_k(T)-V_{k-1}(T)
    \\
    &=
    -c_1|k|+c_1|k-1|
    \\
    &=
    c_1.
\end{align}
At $k=0$, we similarly have
\[
\Delta V_0(T)
=
V_0(T)-V_{-1}(T)
=
c_1.
\]
Therefore,
\begin{equation}
    \Delta V_k(T)=c_1,
    \qquad
    k\leq0.
    \label{eq:linear_terminal_marginal_value}
\end{equation}

Hence, each additional unresolved oversold ticket carries the same marginal
terminal cost.

% ------------------------------------------------------------
\subsubsection{Quadratic penalty}
% ------------------------------------------------------------

For $c_2>0$, define
\begin{equation}
    \Psi_{\mathrm{quad}}(k)
    =
    c_2|k|^2,
    \qquad
    k<0.
    \label{eq:quadratic_penalty}
\end{equation}

For $k<0$, the terminal condition gives
$V_k(T)=-c_2|k|^2$, so that
\begin{align}
    \Delta V_k(T)
    &=
    V_k(T)-V_{k-1}(T)
    \\
    &=
    -c_2|k|^2+c_2|k-1|^2
    \\
    &=
    c_2(2|k|+1).
\end{align}
At $k=0$,
\[
\Delta V_0(T)
=
V_0(T)-V_{-1}(T)
=
c_2,
\]
which is consistent with the same expression. Therefore,
\begin{equation}
    \Delta V_k(T)
    =
    c_2(2|k|+1),
    \qquad
    k\leq0.
    \label{eq:quadratic_terminal_marginal_value}
\end{equation}

Unlike the linear specification, the quadratic penalty generates a marginal
terminal cost that increases with the depth of overbooking. This distinction is
central to the numerical comparison of the two penalty structures.

% ============================================================
\subsection{Extended Kolmogorov equation}
\label{sec:extended-kfe}
% ============================================================

To close the MFG system, the optimal individual policy must be consistent with
the evolution of the aggregate inventory distribution. Let $m_k(t)$ denote the
proportion of firms at state
\[
k\in\mathcal{X}_{\mathrm{ext}}.
\]

For an interior state
\[
k\in\{-M+1,\dots,K-1\},
\]
the forward equation is
\begin{equation}
    \frac{dm_k}{dt}
    =
    \lambda_{k+1}^*(t)m_{k+1}(t)
    -
    \lambda_k^*(t)m_k(t)
    +
    \mu(k-1,t)m_{k-1}(t)
    -
    \mu(k,t)m_k(t),
    \label{eq:extended_kfe}
\end{equation}
where
\begin{equation}
    \lambda_k^*(t)
    =
    \lambda
    \left(
        p_k^*(t),
        \bar{p}(t),
        m(t)
    \right)
\end{equation}
is the optimal sales intensity.

The first two terms describe downward sales transitions, while the final two
terms describe upward cancellation transitions.

% ------------------------------------------------------------
\subsubsection{Boundary dynamics and mass conservation}
\label{sec:extended-mass}
% ------------------------------------------------------------

At the lower boundary $k=-M$, further sales are suppressed while cancellations
remain possible:
\begin{equation}
    \frac{dm_{-M}}{dt}
    =
    \lambda_{-M+1}^*(t)m_{-M+1}(t)
    -
    \mu(-M,t)m_{-M}(t).
    \label{eq:kfe_lower_boundary}
\end{equation}

At the upper boundary $k=K$,
\begin{equation}
    \frac{dm_K}{dt}
    =
    -\lambda_K^*(t)m_K(t)
    +
    \mu(K-1,t)m_{K-1}(t).
    \label{eq:kfe_upper_boundary}
\end{equation}

Together, Equations~\eqref{eq:extended_kfe}--\eqref{eq:kfe_upper_boundary}
conserve total probability mass:
\begin{equation}
    \frac{d}{dt}
    \sum_{k=-M}^{K}
    m_k(t)
    =
    0.
\end{equation}
Hence, an initial distribution satisfying
\[
\sum_{k=-M}^{K}m_k(0)=1
\]
remains normalized for all $t\in[0,T]$.

Firms in states $k=-M+1,\dots,K$ remain sales-active, whereas further sales are
suppressed at $k=-M$. We therefore define
\begin{equation}
    \eta_{\mathrm{ext}}(t)
    =
    \sum_{k=-M+1}^{K}m_k(t)
    =
    1-m_{-M}(t).
    \label{eq:extended_active_mass}
\end{equation}

The corresponding mean quoted price among sales-active firms is
\begin{equation}
    \bar{p}(t)
    =
    \frac{
        \displaystyle
        \sum_{k=-M+1}^{K}
        m_k(t)p_k^*(t)
    }{
        \eta_{\mathrm{ext}}(t)
    },
    \qquad
    \eta_{\mathrm{ext}}(t)>0.
\label{eq:extended_mean_price}
\end{equation}
The normalization by $\eta_{\mathrm{ext}}(t)$ conditions on sales-active firms,
preventing changes in active mass alone from mechanically changing the mean
quoted price.

Together with the extended HJB equation, these forward dynamics define the
coupled HJB--Kolmogorov system.

\begin{definition}[Selected mean field equilibrium]
\label{def:selected-mfg-equilibrium}
A selected Markovian mean field equilibrium is a quadruple
$(V,\alpha^*,m,\bar p)$ such that:
\begin{enumerate}
 \item[(i)] for the prescribed flow $(m,\bar p)$, the feedback $\alpha^*$ is
 optimal for the representative-firm problem and $V$ is its value function;
 \item[(ii)] $m(0)=m^0$, and $m$ is the law of the optimally controlled
 inventory process, equivalently, it solves the forward Kolmogorov equations
 with the rates induced by $\alpha^*$; and
 \item[(iii)] whenever $\eta_{\mathrm{ext}}(t)>0$, the mean-quoted-price
 consistency condition holds:
 \begin{equation}
  \bar p(t)=\frac{1}{\eta_{\mathrm{ext}}(t)}
  \sum_{k=-M+1}^{K}m_k(t)\alpha^*(t,k),
  \qquad
  \eta_{\mathrm{ext}}(t)=\sum_{k=-M+1}^{K}m_k(t).
  \label{eq:equilibrium_price_consistency}
 \end{equation}
\end{enumerate}
The qualifier ``selected'' refers to the tie-breaking rule in
Equation~\eqref{eq:selected_price} for the zero-demand region.
\end{definition}

This is a non-atomic MFG equilibrium: a representative firm cannot improve by
a unilateral deviation when the population flow is held fixed. It is not an
exact Nash equilibrium of a finite-player game. Establishing an
$\varepsilon_N$-Nash approximation would require a separately specified
$N$-firm model and is not claimed here.

% ============================================================
\subsection{Structural properties}
\label{sec:theoretical-insights}
% ============================================================

The terminal penalty provides a useful interpretation of pricing behavior near
and within the negative inventory region. Recall that
\[
\Delta V_k(t)
=
V_k(t)-V_{k-1}(t)
\]
is the marginal value of inventory.

For both penalty specifications,
\[
\Delta V_k(T)>0,
\qquad k\leq0,
\]
while under the quadratic specification the terminal marginal value increases
with the depth of overbooking. Under the adopted selection, a larger marginal
value raises the quoted price. The selection-independent implication is that
the unique optimal sales intensity in Equation~\eqref{eq:unique_sales_intensity}
decreases and may become zero, thereby increasing the opportunity cost of an
additional sale. Through the backward HJB equation, this terminal effect
propagates to earlier times.

A sufficiently steep convex penalty can therefore discourage progressively
deeper overselling without introducing an additional hard boundary inside the
admissible state space. This mechanism motivates the endogenous soft-boundary
behavior examined numerically in Section~\ref{sec:penalty-intensity}.

% ============================================================
\subsection{Existence and equilibrium selection}
\label{sec:extended-equilibrium}
% ============================================================

The compact-action assumption is standard in finite-state MFG existence
theory. Interaction through the distribution of controls requires additional
care; see \citet{cecchin2020probabilistic} for finite-state probabilistic
methods, \citet{carmona2021extended} for extended finite-state MFGs with
control interaction, and \citet{garcia2025control} for controlled-intensity
MFGs. Because the present model also contains bidirectional inventory jumps
and an active-mass normalization, we give a specialized fixed-point argument
in Appendix~\ref{app:existence} rather than claim a direct application of any
one of these results.

\begin{proposition}[Existence of a selected mean field equilibrium]
\label{prop:selected-equilibrium-existence}
Fix a finite time horizon $T>0$. Suppose that $\mathcal X_{\mathrm{ext}}$ is finite,
$\mathcal A=[0,P_{\max}]$ with $P_{\max}\geq1$, the cancellation rates
$\mu(k,\cdot)$ are bounded and continuous on $[0,T]$, the terminal penalty is
finite on $\mathcal X_{\mathrm{ext}}$, and the competition parameter satisfies
$\epsilon\geq0$. Let the initial distribution
$m^0\in\mathcal P(\mathcal X_{\mathrm{ext}})$ satisfy
$\eta_{\mathrm{ext}}(0)=\sum_{k=-M+1}^{K}m_k^0>0$. Under the feedback selection
in Equation~\eqref{eq:selected_price}, the extended model admits at least one
selected Markovian mean field equilibrium on $[0,T]$.
\end{proposition}

For a prescribed mean-field flow, the finite-state HJB system has a unique
value function. The induced optimal sales intensity is also unique, even where
the quoted price is not. For fixed $(V,m)$, the contemporaneous
mean-quoted-price consistency equation is a contraction and therefore has a
unique solution.

These componentwise statements do not imply uniqueness of the complete mean
field equilibrium. The baseline model proves uniqueness under a
small-interaction condition, while general finite-state MFG results commonly
use restrictions such as a short horizon or monotonicity. In the present
extension, terminal penalties can produce zero optimal sales intensities, so
the strictly positive-intensity argument used in the baseline model does not
apply directly. We therefore do not claim global uniqueness for arbitrary
interaction strength, horizon, or quoted-price selection. The numerical
fixed-point convergence reported below is a computational observation, not a
proof of continuous-time equilibrium uniqueness.

\section{Numerical Results}
\label{sec:results}
% ============================================================

This section examines how terminal penalties, cancellation dynamics, and
mean-field coupling affect equilibrium pricing and inventory behavior in the
extended model.

All quoted prices in this section are generated by the tie-breaking selection
in Equation~\eqref{eq:selected_price}, with $P_{\max}=12$. The largest value
over all reported runs and terminal grid nodes is $10.563255$, attained at
$(k,t)=(-99,200)$ in the quadratic-penalty stress case with $M=100$; hence the
cap is nonbinding.

The coupled HJB--Kolmogorov system is solved by relaxed Picard iteration
(relaxation parameter $0.5$) and a four-stage Runge--Kutta scheme. Prescribed
paths or rates are linearly interpolated and evaluated at each Runge--Kutta
stage. Unless stated otherwise, $\Delta t=0.1$. For the unrelaxed update
$(\widehat{\bar p}^{(n+1)},\widehat m^{(n+1)})$, iteration stops when both
$\lVert\widehat{\bar p}^{(n+1)}-\bar p^{(n)}\rVert_\infty<10^{-7}$ and
$\max_j\lVert\widehat m^{(n+1)}(t_j)-m^{(n)}(t_j)\rVert_1<10^{-7}$.

The benchmark parameters are
\[
K=100,\quad T=200,\quad r=0.04,\quad M=20,\quad
\epsilon=0.4,\quad \gamma=0.005,\quad R=0.3.
\]

The benchmark initial inventory distribution assigns probability $0.1$ to
each state in $\{20,21,22,23,24,50,51,52,53,54\}$ and zero elsewhere; the
alternative assigns equal probability to every nonnegative state.

The parameter $M$ determines the lower truncation boundary of the finite
inventory state space. All benchmark and graphical experiments use $M=20$.
For the low-cancellation terminal-risk comparison reported in
Section~\ref{sec:quantitative-metrics}, the state space is enlarged to $M=80$
to reduce truncation effects in the lower tail. The sensitivity of these
terminal-risk measures to the choice of $M$ is examined separately in
Section~\ref{sec:numerical-robustness}.

For the terminal penalties, we use
\[
c_1=2.0
\qquad\text{and}\qquad
c_2=0.1
\]
for the linear and quadratic specifications, respectively. The coefficients are
chosen so that the two penalty functions coincide at the common reference state
$k=-20$:
\begin{equation}
    \Psi_{\mathrm{lin}}(-20)
    =
    \Psi_{\mathrm{quad}}(-20)
    =
    40.
    \label{eq:matched_penalty_boundary}
\end{equation}
This matching provides a common reference level while allowing the numerical
comparison to focus on differences in penalty shape. All parameters are treated
as normalized benchmark inputs rather than empirical estimates.
%============================================================
\subsection{Terminal penalties and marginal risk}
\label{sec:marginal-risk}
% ============================================================

The marginal opportunity cost of inventory,
\[
\Delta V_k(t)
=
V_k(t)-V_{k-1}(t),
\]
plays a central role in the optimal pricing rule. We begin by examining its
terminal value at $t=T$, where the overbooking penalty is realized.

\begin{figure}[htbp]
    \centering
    \includegraphics[width=0.5\textwidth]{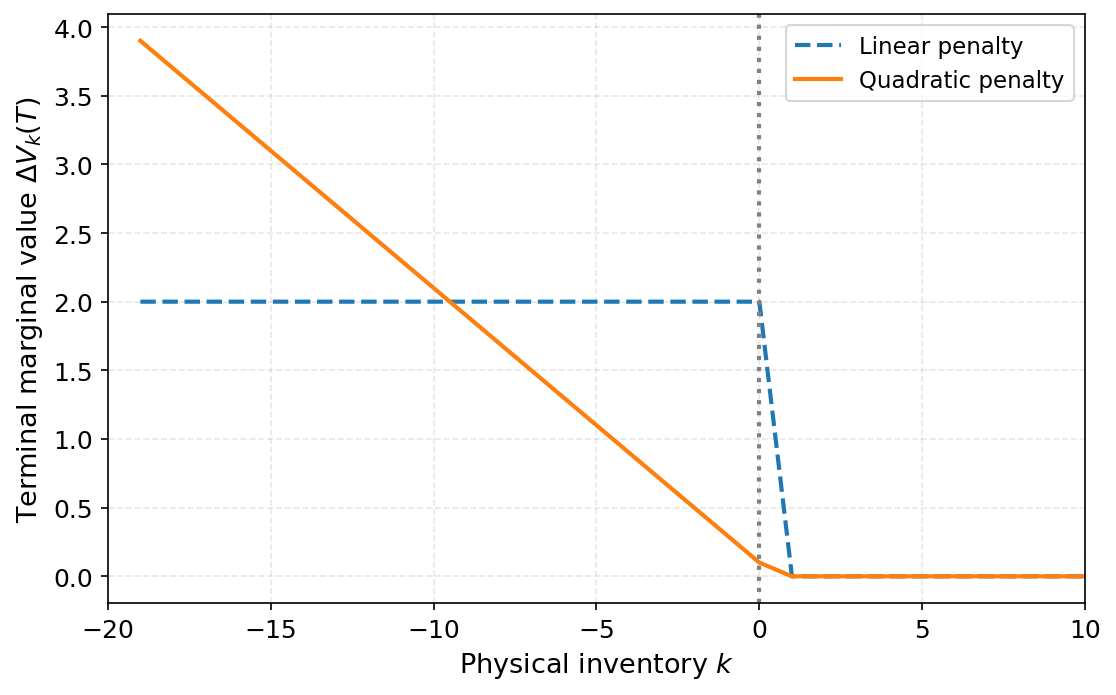}
    \caption{Terminal marginal value $\Delta V_k(T)$ under the linear and quadratic terminal penalty specifications}
    \label{fig:penalty_cliff}
\end{figure}

For positive inventory states, the absence of terminal salvage value gives
\[
\Delta V_k(T)=0,
\qquad
k>0.
\]
At $k=0$, however, an additional sale moves the firm into the negative
inventory region and therefore incurs the marginal terminal overbooking cost.

Under the linear penalty, this marginal cost is constant:
\[
\Delta V_k(T)=c_1,
\qquad
k\leq0.
\]
Thus, each additional unresolved oversold ticket carries the same marginal
terminal cost regardless of the current depth of overbooking.

Under the quadratic penalty, the marginal cost increases as inventory becomes
more negative:
\[
\Delta V_k(T)
=
c_2(2|k|+1),
\qquad
k\leq0.
\]

As shown in Figure~\ref{fig:penalty_cliff}, the linear penalty therefore
generates a constant marginal terminal cost in the oversold region, whereas
the quadratic penalty places progressively greater marginal cost on deeper
negative inventory states. This difference provides the main mechanism through
which convex terminal penalties can discourage severe overbooking and affects
the near-departure pricing behavior examined next.

% ============================================================
\subsection{Near-departure pricing}
\label{sec:near-departure-pricing}
% ============================================================

The differences in terminal marginal value propagate backward through the
extended HJB equation and influence the selected pricing strategy before
departure. We examine the pricing policy at $t=190$, shortly before the
terminal time $T=200$.

\begin{figure}[htbp]
    \centering
    \includegraphics[width=0.5\textwidth]{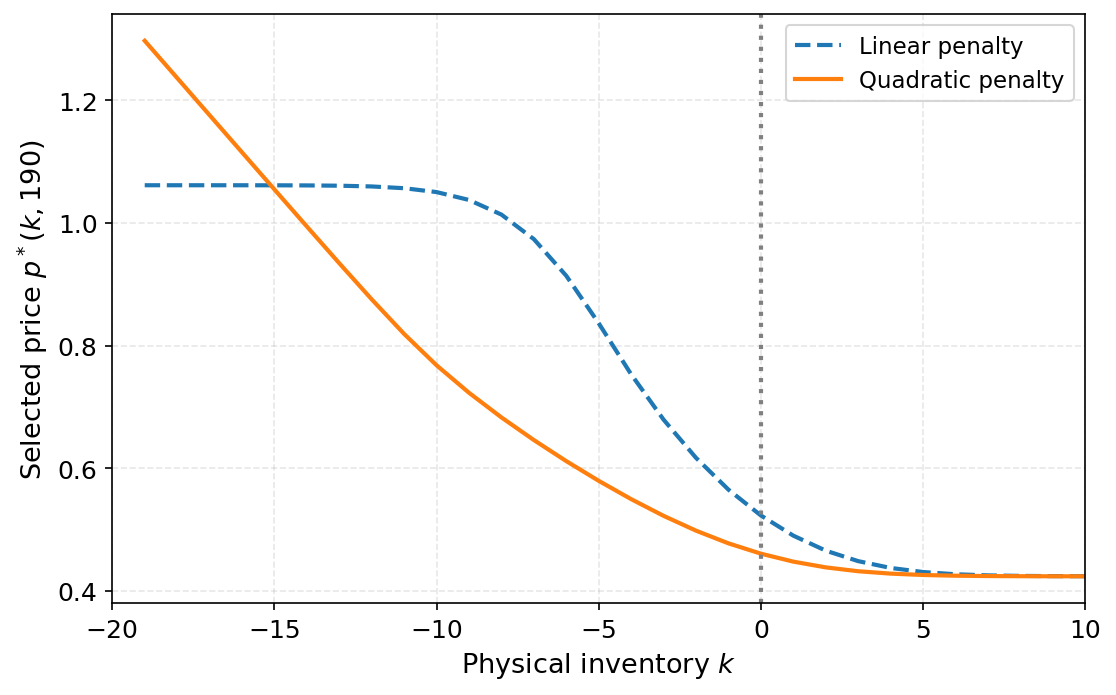}
    \caption{Selected optimal prices at $t=190$ under the linear and quadratic terminal penalty specifications}
    \label{fig:penalty_premium}
\end{figure}

Figure~\ref{fig:penalty_premium} shows that for positive inventory states, the
two penalty regimes generate similar pricing behavior. Although firms in these
states still have remaining physical capacity, terminal overbooking penalties
can affect current pricing through continuation values, since subsequent sales
may move inventory into the negative region.

The pricing behavior changes once inventory enters the negative region. Under
the linear penalty, the constant marginal terminal cost generates a relatively
high price over a range of mildly oversold states. Under the quadratic penalty,
the price increases more gradually near the boundary of the negative region,
because the marginal terminal cost is still relatively small when overselling
is shallow.

As inventory becomes more negative, however, the quadratic marginal terminal
cost rises with the depth of overbooking. Under the adopted selection, the
quadratic penalty therefore induces a higher quoted price than the linear
penalty in sufficiently deep negative states. The selection-independent
implication is that the unique optimal sales intensity decreases and may become
zero, discouraging further movement into the negative tail. This intensity
response is the soft-boundary mechanism generated by the convex terminal
penalty.

% ============================================================
\subsection{Sensitivity to cancellation rates}
\label{sec:cancellation-sensitivity}
% ============================================================

We next examine how the equilibrium responds to changes in the cancellation
rate while maintaining the quadratic terminal penalty. The cancellation rate
controls the frequency of upward inventory transitions and therefore the
ability of firms to recover inventory before departure.

Because the strongest effects arise in the low-cancellation region, we compare
\[
\gamma
\in
\{0.001,0.002,0.005\},
\]
where $\gamma=0.005$ is the benchmark specification. All other parameters,
including the lower inventory boundary $M=20$, are kept at their benchmark
values.

To summarize the aggregate state of the market, we report the sales-active
mass
\begin{equation}
    \eta_{\mathrm{ext}}(t)
    =
    \sum_{k=-M+1}^{K} m_k(t)
    =
    1-m_{-M}(t),
\end{equation}
which measures the fraction of firms whose current inventory lies above the
lower overbooking boundary and for which further sales are admissible.

\begin{figure}[htbp]
    \centering
    \includegraphics[width=0.95\textwidth]{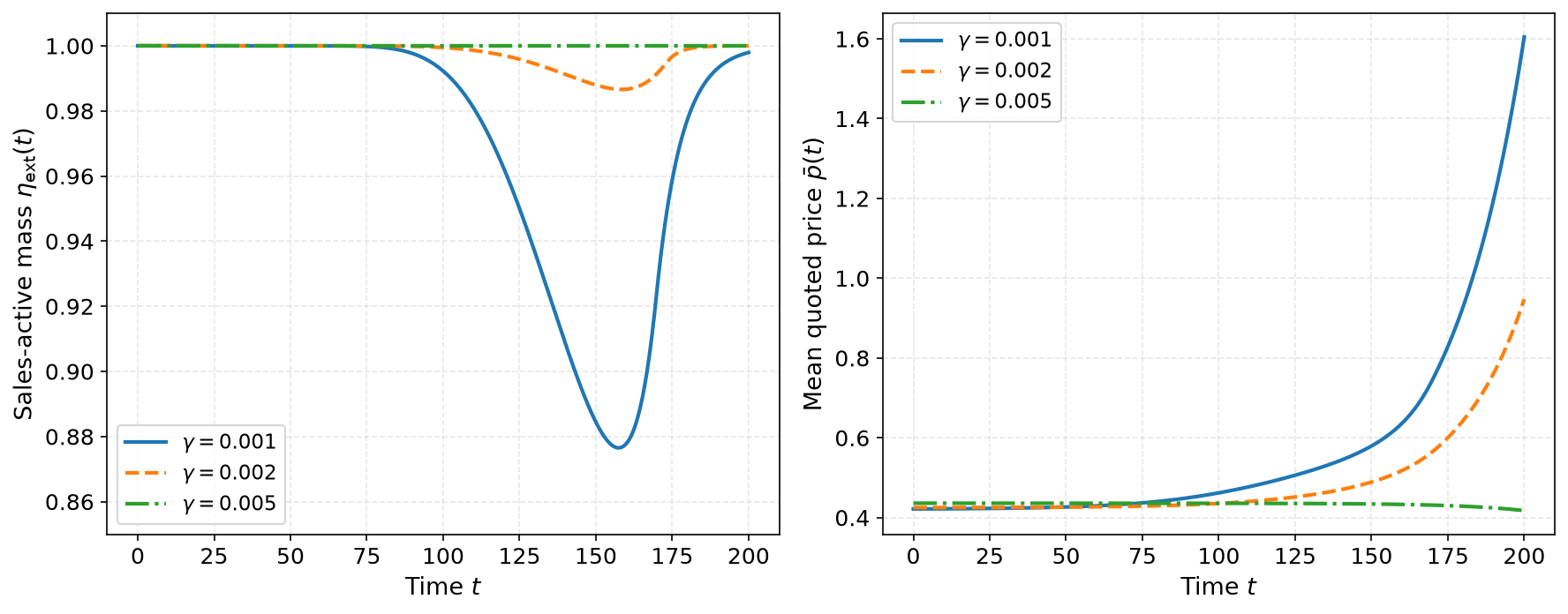}
    \caption{Sensitivity of sales-active mass and mean quoted prices to the cancellation rate $\gamma$}
    \label{fig:cancellation_sensitivity}
\end{figure}

Figure~\ref{fig:cancellation_sensitivity} shows a pronounced response to the
cancellation rate. When $\gamma=0.001$, cancellation-driven inventory recovery
is weak, so a substantial fraction of firms is pushed toward the lower
overbooking boundary during the selling horizon. As a result, the sales-active
mass falls markedly, reaching a minimum of about $0.877$. Increasing the
cancellation rate to $\gamma=0.002$ greatly mitigates this effect, with the
minimum sales-active mass remaining near $0.987$, while under the benchmark
value $\gamma=0.005$ the sales-active mass stays essentially equal to one
throughout the horizon.

The pricing response moves in the opposite direction. Lower cancellation rates
make inventory recovery less likely and increase the effective risk of
remaining in a negative inventory state near departure. Firms therefore raise
prices more strongly in the low-cancellation cases in order to slow further
sales. This effect is particularly strong when $\gamma=0.001$, for which the
mean quoted price rises substantially as the terminal date approaches. Under
$\gamma=0.002$, the increase is more moderate, and under the benchmark value
$\gamma=0.005$ the mean quoted price remains comparatively stable over time.

These results show that the cancellation rate is an important determinant of
both aggregate market activity and defensive pricing in the extended model.
Low cancellation rates make recovery from overselling more difficult, pushing
more firms toward the lower boundary and inducing a stronger increase in
equilibrium prices. Higher cancellation rates provide a stronger recovery
mechanism and thereby reduce both boundary pressure and the need for aggressive
near-terminal price adjustments. The quantitative magnitude of these effects
remains specific to the chosen calibration.
% ============================================================
\subsection{Time-varying cancellations}
\label{sec:anticipatory-overbooking}
% ============================================================

The constant cancellation rate provides a useful benchmark, but the timing of
cancellations may also affect current pricing and inventory behavior. We
therefore consider the time-dependent specification introduced in
Equation~\eqref{eq:gamma_time_dependent}.

For this experiment, we set
\[
\gamma_{\mathrm{base}}=0.001,
\qquad
\gamma_{\mathrm{surge}}=0.02,
\qquad
\beta=4,
\]
so that
\[
\gamma(t)
=
0.001
+
0.02
\left(
    \frac{t}{T}
\right)^4.
\]
Cancellation activity is therefore relatively limited early in the selling
horizon and increases toward departure.

The time-average cancellation rate is
\begin{equation}
    \frac{1}{T}
    \int_0^T \gamma(t)\,dt
    =
    \gamma_{\mathrm{base}}
    +
    \frac{\gamma_{\mathrm{surge}}}{\beta+1}
    =
    0.005,
    \label{eq:average_dynamic_gamma}
\end{equation}
which coincides with the constant-rate benchmark $\gamma=0.005$. The comparison
therefore examines the effect of cancellation timing while keeping the
time-average per-ticket cancellation rate unchanged. The expected total number
of cancellations need not be the same, since the aggregate cancellation
intensity also depends on the endogenous number of outstanding bookings.

\begin{figure}[htbp]
    \centering
    \includegraphics[width=0.95\textwidth]{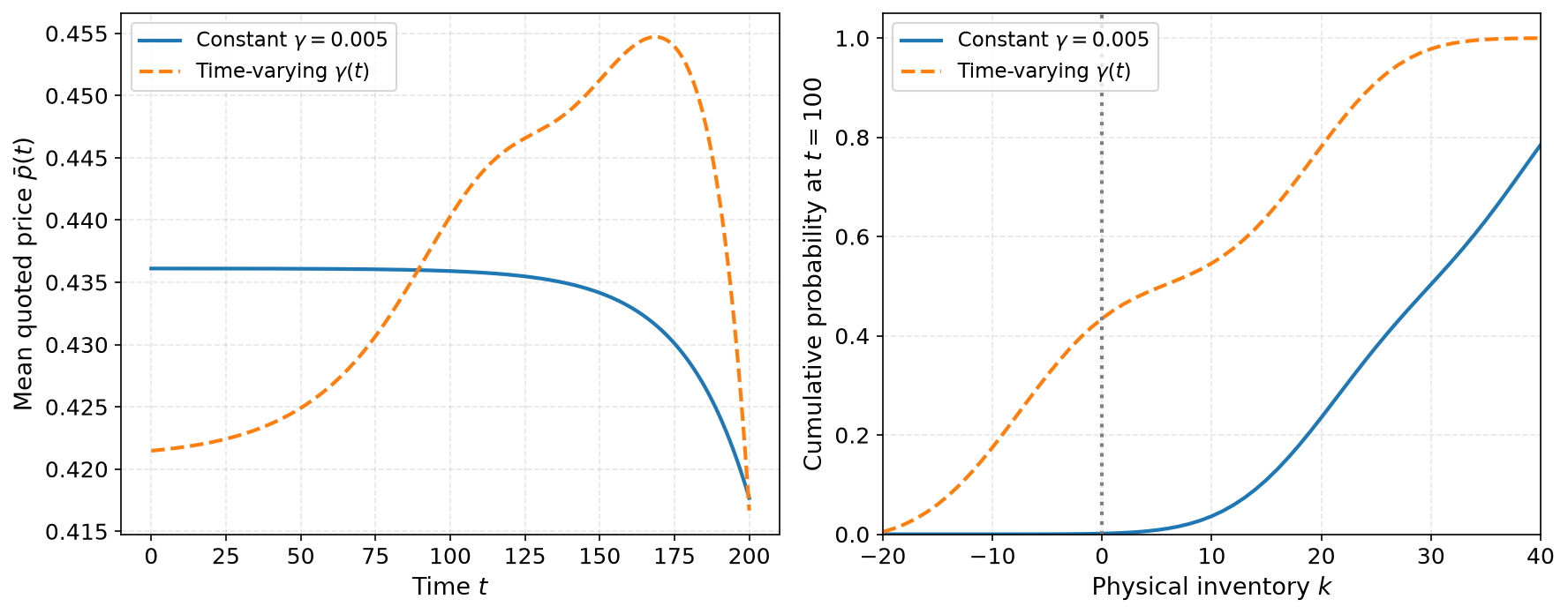}
    \caption{Mean quoted prices and inventory distributions under constant and time-varying cancellation rates, with the right panel evaluated at $t=100$}
    \label{fig:dynamic_results}
\end{figure}

Figure~\ref{fig:dynamic_results} shows that the timing of cancellation
opportunities materially changes equilibrium behavior. Early in the selling
horizon, the time-varying specification produces a lower mean quoted price than the
constant-rate benchmark. Although current cancellation activity is relatively
low, firms anticipate stronger inventory recovery later in the horizon, which
reduces the effective cost of carrying temporary overbooking exposure.

As the horizon progresses, however, the mean quoted price under the time-varying
specification rises above the constant-rate benchmark. By this stage, the lower
cancellation activity earlier in the horizon has allowed inventory to deplete
more rapidly, increasing the need for defensive pricing as departure
approaches.

The inventory distributions at $t=100$ are consistent with this interpretation.
The cumulative distribution under the time-varying specification is shifted
toward lower inventory states, indicating substantially greater inventory
depletion and a larger mass near or below the physical capacity boundary than
under the constant-rate benchmark.

The comparison therefore shows that cancellation timing matters even when the
average cancellation rate is unchanged. Anticipated future recovery can support
greater temporary inventory depletion early in the horizon, while the resulting
inventory position can subsequently induce stronger near-terminal pricing
responses. This experiment should be interpreted as a timing comparison rather
than as an empirical description of airline cancellation behavior.

% ============================================================
\subsection{Penalty intensity and effective boundary}
\label{sec:penalty-intensity}
% ============================================================

We next examine how the strength of the quadratic terminal penalty affects
pricing behavior in the negative inventory region. Rather than establishing a
formal limiting result, the purpose of this experiment is to investigate
numerically how increasing penalty intensity changes the incentive to enter
deeper oversold states.

We evaluate the selected pricing policy at
\[
t=190
\]
for
\[
c_2
\in
\{0.01,0.05,0.1,0.2\}.
\]

\begin{figure}[htbp]
    \centering
    \includegraphics[width=0.7\textwidth]{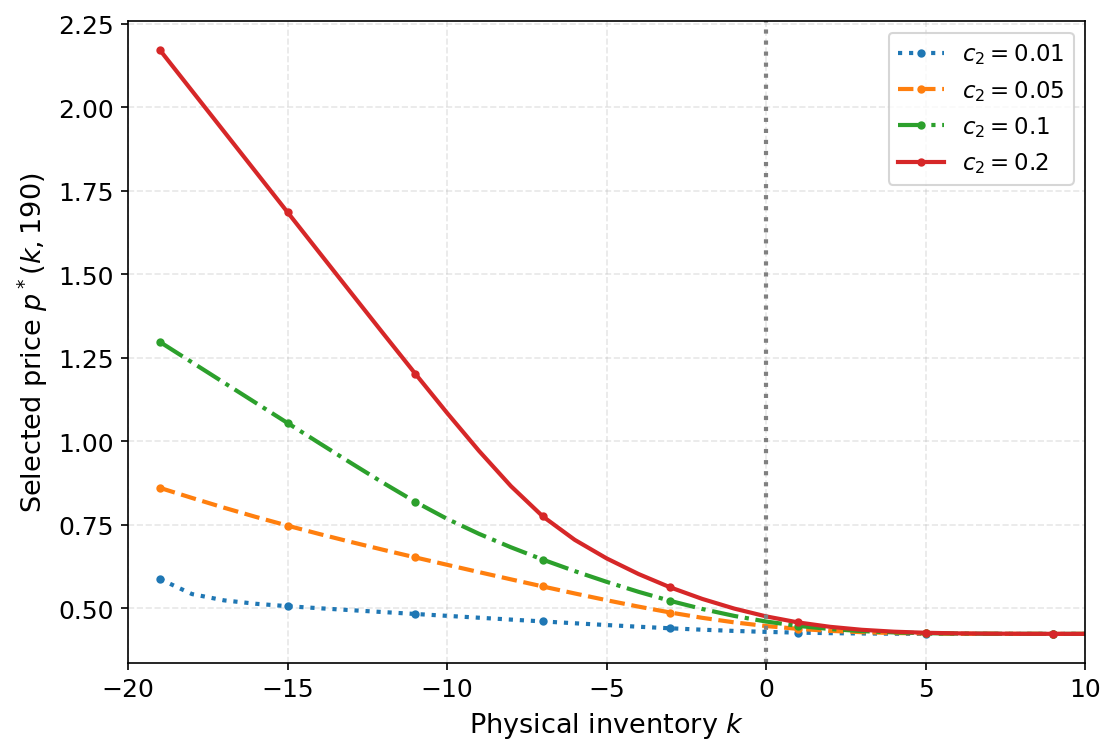}
    \caption{Selected optimal prices at $t=190$ for different quadratic penalty coefficients $c_2$}
    \label{fig:penalty_intensity}
\end{figure}

Figure~\ref{fig:penalty_intensity} shows that the effect of the penalty
coefficient is concentrated in the negative inventory region. For nonnegative
inventory states, the pricing policies remain close across the different
values of $c_2$. Once inventory becomes negative, however, the policies begin
to separate, and the difference becomes increasingly pronounced as the depth
of overbooking grows.

For a small penalty coefficient, such as $c_2=0.01$, the selected price rises
only moderately as inventory becomes more negative. Increasing $c_2$ produces
a progressively steeper response under the adopted selection. In particular,
deeply oversold states are associated with substantially higher selected
prices when the quadratic penalty is stronger.

The robust implication is the associated reduction in the uniquely determined
sales intensity. The quadratic penalty therefore does not impose a new hard
capacity constraint, but increasingly discourages movement into deeper
oversold states as $c_2$ rises. This numerical mechanism motivates the
interpretation of the convex terminal penalty as generating an endogenous soft
boundary; it is not asserted as a general limiting theorem.

% ============================================================
\subsection{Quantitative metrics and robustness}
\label{sec:quant-robustness}
% ============================================================

We complement the preceding graphical analysis with quantitative metrics and
numerical robustness checks. These results are intended to compare structural
properties of the model rather than to provide an empirical calibration to
specific airline routes.

% ------------------------------------------------------------
\subsubsection{Quantitative metrics}
\label{sec:quantitative-metrics}
% ------------------------------------------------------------

We compare the three penalty regimes under the low-cancellation terminal-risk
stress case
\[
\gamma=0.001,
\qquad
\epsilon=0.4.
\]
Because this specification can generate substantially deeper negative inventory
states than the benchmark case, the lower truncation boundary is enlarged to
$M=80$ for this comparison.

To distinguish the prevalence, conditional severity, and lower-tail risk of
overselling, we report three terminal-distribution measures. The terminal
oversold mass is
\begin{equation}
    O(T)
    =
    \sum_{k=-M}^{-1} m_k(T)
    =
    \mathbb{P}(X_T<0).
    \label{eq:terminal_oversold_mass}
\end{equation}

Conditional on overselling, the expected depth is
\begin{equation}
    D_{\mathrm{cond}}(T)
    =
    \frac{
        \sum_{k=-M}^{-1}(-k)m_k(T)
    }{
        \sum_{k=-M}^{-1}m_k(T)
    }
    =
    \mathbb{E}
    \left[
        -X_T \mid X_T<0
    \right],
    \label{eq:conditional_overselling_depth}
\end{equation}
whenever $O(T)>0$. The $5\%$ lower quantile is defined by
\begin{equation}
    q_{0.05}(T)
    =
    \min
    \left\{
        k\in\mathcal{X}_{\mathrm{ext}}
        :
        \sum_{j=-M}^{k}m_j(T)\geq0.05
    \right\}.
    \label{eq:terminal_lower_quantile}
\end{equation}

A larger value of $D_{\mathrm{cond}}(T)$ indicates deeper overselling among
oversold firms, while a more negative value of $q_{0.05}(T)$ indicates a more
severe lower tail.

The expected net value reported below is the initial model value averaged over
the initial inventory distribution:
\begin{equation}
    \mathcal V_0=\sum_{k=-M}^{K}m_k(0)V_k(0).
    \label{eq:expected_net_value}
\end{equation}
It includes expected discounted ticket revenues, cancellation-related costs,
and the terminal overbooking penalty.

\begin{table}[htbp]
    \centering
    \small
    \renewcommand{\arraystretch}{1.15}
    \setlength{\tabcolsep}{3pt}

    \begin{tabular}{@{}lccccc@{}}
        \toprule

        \textbf{Penalty Regime}
        &
        \makecell{\textbf{Mean}\\\textbf{Quoted}\\\textbf{Price}\\$\bar p(190)$}
        &
        \makecell{\textbf{Expected}\\\textbf{Net Value}}
        &
        \makecell{\textbf{Terminal}\\\textbf{Oversold Mass}}
        &
        \makecell{\textbf{Conditional}\\\textbf{Depth}}
        &
        \makecell{\textbf{5\% Terminal}\\\textbf{Quantile}}
        \\

        \midrule

        No penalty
        & 0.418
        & 3.815
        & 0.969
        & 27.585
        & $-51$
        \\

        Linear ($c_1=2.0$)
        & 1.091
        & 3.809
        & 0.765
        & 16.022
        & $-33$
        \\

        Quadratic ($c_2=0.1$)
        & 1.326
        & 3.812
        & 0.967
        & 11.689
        & $-24$
        \\

        \bottomrule
    \end{tabular}

    \caption{Quantitative comparison under the low-cancellation specification
    $\gamma=0.001$ and $\epsilon=0.4$, using $M=80$}
    \label{tab:metrics}
\end{table}
Table~\ref{tab:metrics} shows that the penalty specifications affect different
dimensions of overbooking. Without a terminal penalty, $96.9\%$ of firms end
the horizon in negative inventory, with a conditional overselling depth of
$27.585$ and a $5\%$ terminal quantile of $-51$.

The linear penalty most strongly reduces the prevalence of terminal
overselling, lowering the oversold mass to $0.765$. The quadratic penalty
instead produces the smallest conditional depth, $11.689$, and the least
negative $5\%$ terminal quantile, $-24$, despite leaving the overall oversold
mass relatively high at $0.967$. Thus, under the present calibration, the
linear penalty mainly reduces the probability of ending in negative inventory,
whereas the quadratic penalty more strongly limits the severity of overselling
once it occurs. This pattern is consistent with the soft-boundary mechanism
discussed in Section~\ref{sec:penalty-intensity}.

The expected net values remain close across the three regimes, ranging from
$3.809$ to $3.815$. This similarity is specific to the present calibration and
should not be interpreted as a general equivalence between the penalty
specifications.

\subsubsection{Numerical robustness}
\label{sec:numerical-robustness}
% ------------------------------------------------------------

The $M=20$ boundary materially restricts the stress-case negative tail, so
Table~\ref{tab:metrics} uses $M=80$. Increasing it to $100$ changes
$\bar p(190)$, expected net value, $O(T)$, and $D_{\mathrm{cond}}(T)$ by less
than $0.0004\%$; the three $5\%$ quantiles remain $-51,-33,-24$. The largest
lower-boundary mass over time is below $5.33\times10^{-6}$ for $M=80$ and
$3.28\times10^{-10}$ for $M=100$. Likewise, for
$\Delta t\in\{0.2,0.1,0.05\}$, refining $0.1$ to $0.05$ changes those four
continuous metrics by less than $0.0004\%$ and leaves the quantiles unchanged.

For each of the fourteen principal cases, Picard iteration from the constant
paths $\bar p^{(0)}\equiv0,0.5,12$ gives maximum pairwise differences
$5.10\times10^{-8}$ in mean quoted price and $4.24\times10^{-8}$ in the
distribution path (maximum-time $\ell^1$ norm). This supports numerical
stability for the tested initializations but does not prove equilibrium
uniqueness. Finally, without clipping or renormalization, the maximum mass
error over all verified runs and Picard iterates is $2.23\times10^{-15}$;
the minimum computed probability is zero, so no negative values occur.

% ------------------------------------------------------------
\subsubsection{Robustness to the initial distribution}
\label{sec:initial-distribution-robustness}
% ------------------------------------------------------------

We compare the benchmark bimodal initial distribution with a uniform
distribution over the nonnegative states, keeping $\epsilon=0.4$,
$\gamma=0.005$, and all other parameters fixed.

\begin{figure}[htbp]
    \centering
    \includegraphics[width=0.95\textwidth]{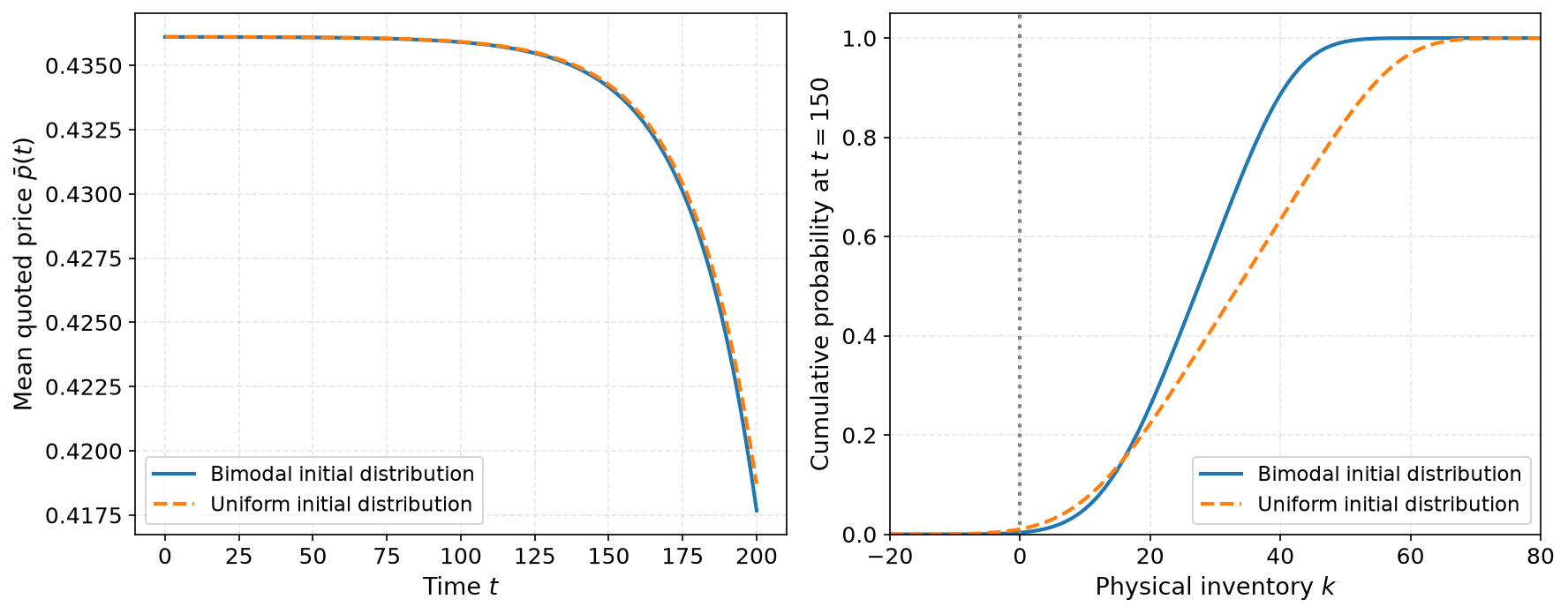}
    \caption{Mean quoted prices and inventory distributions under bimodal and
    uniform initial conditions, with the right panel evaluated at $t=150$}
    \label{fig:initial_distribution_robustness}
\end{figure}

Figure~\ref{fig:initial_distribution_robustness} shows close mean-quoted-price
paths with the same terminal adjustment. The inventory CDFs at $t=150$ cross:
the uniform case is lower over much of the range but higher at some
low-inventory thresholds, so it is not a uniform rightward shift. Thus,
aggregate pricing is relatively robust here, whereas the distributional shape
remains sensitive to the initial condition.
% ============================================================
\subsection{Mean-field coupling and overbooking incentives}
\label{sec:competition-impact}
% ============================================================

We next examine how the strength of mean-field coupling affects equilibrium
pricing and inventory depletion under the quadratic terminal penalty. We
compare
\[
\epsilon
\in
\{0,0.2,0.4\},
\]
where $\epsilon=0$ provides an uncoupled benchmark and larger values correspond
to stronger mean-field interaction.

As discussed in Section~\ref{sec:demand-competition}, changing $\epsilon$
affects both the weight placed on the mean quoted price and the effective
demand component of the sales intensity. The comparison should therefore be
interpreted as a model-specific comparative static rather than as a general
statement about market structure.

\begin{figure}[htbp]
    \centering
    \includegraphics[width=0.95\textwidth]{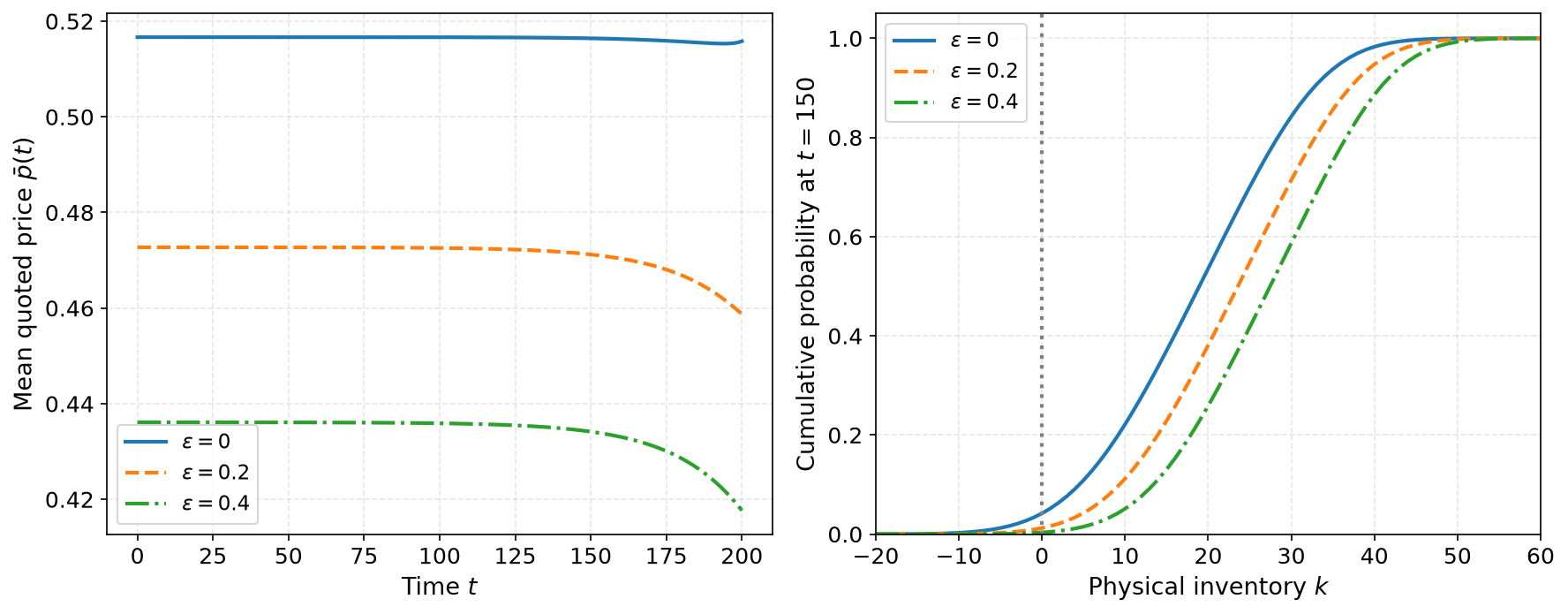}
    \caption{Mean quoted prices and inventory distributions for
    $\epsilon\in\{0,0.2,0.4\}$, with the right panel evaluated at $t=150$}
    \label{fig:competition_impact}
\end{figure}

Figure~\ref{fig:competition_impact} shows a clear ordering across the three
coupling levels. Under the adopted price selection, the mean quoted price
decreases as $\epsilon$ increases.
At $t=150$, the corresponding values are approximately $0.516$, $0.471$, and
$0.434$ for $\epsilon=0$, $0.2$, and $0.4$, respectively.

The inventory distributions show a corresponding shift toward higher inventory
states as $\epsilon$ increases, indicating slower inventory depletion. At
$t=150$, the fraction of firms in negative inventory decreases from
approximately $0.0336$ for $\epsilon=0$ to $0.0093$ for $\epsilon=0.2$ and
$0.0023$ for $\epsilon=0.4$.

Thus, stronger mean-field coupling is associated with both lower equilibrium
prices and slower inventory depletion under the present specification. This
pattern reflects the demand normalization inherited from the baseline model.
Increasing $\epsilon$ changes not only the firm's exposure to the mean quoted
price but also the effective demand level entering the sales intensity.

The observed ordering should therefore not be interpreted as a general
prediction that stronger competition necessarily reduces overbooking. Rather,
it is a comparative-static result specific to the adopted demand
specification.

\section{Conclusion}
\label{sec:conclusion}
% ============================================================

This paper extends the strict-capacity MFG framework for dynamic ticket pricing
by allowing negative inventory states, stochastic cancellation-driven inventory
recovery, and terminal penalties for unresolved overbooking. The resulting
finite-state model captures the interaction between pricing, temporary
overselling, expected cancellations, and terminal overbooking risk.
We defined the representative-firm problem and the population consistency
conditions and proved existence of a selected Markovian mean field equilibrium
under the stated assumptions, including a compact price set and an explicit
tie-breaking rule. The selected quoted price need not be the unique optimizer
in zero-demand states, although the induced optimal sales intensity is unique.

The numerical results highlight three main mechanisms. First, the quadratic
terminal penalty generates a numerical soft-boundary mechanism: its marginal
cost increases with the depth of overselling, thereby reducing the uniquely
determined optimal sales intensity. Second, the time-varying cancellation
experiment suggests that anticipated future inventory recovery can support
lower early prices and greater temporary overbooking exposure. Third, under
the adopted demand normalization, stronger mean-field coupling is associated
with lower equilibrium prices and slower inventory depletion.
The quantitative results further indicate that linear and quadratic penalties
affect different dimensions of overbooking: the linear penalty more strongly
reduces its prevalence, whereas the quadratic penalty more strongly limits its
depth under the reported calibration.

Several limitations remain. Cancellation intensities are exogenous, the model
contains a single homogeneous seat class, and the parameter values are
stylized rather than empirically calibrated. Future work could incorporate
price-dependent cancellation behavior and multiple fare classes, empirically
calibrate the model, and identify structural or parameter conditions ensuring
uniqueness of the complete mean field equilibrium. The existence result
established here neither requires nor implies such global uniqueness.

\section*{Acknowledgments}

The authors thank Antonis Papapantoleon for valuable comments and suggestions
that helped improve an earlier version of this work.

\appendix
\section{Baseline Model Formulation}
\label{app:baseline}

For completeness, we provide the core mathematical formulation of the baseline Mean Field Game framework introduced by \citet{aydin2026fare}, upon which our extended model is built.

\paragraph{State Dynamics and the HJB Equation.}
In the strict-capacity baseline, the firm's inventory $X_t \in \{0, 1, \dots, K\}$ follows a pure death process driven by a Poisson counting process with intensity $\lambda_t$, such that $dX_t = -\mathbb{I}_{\{X_t>0\}}\,dN_t^\lambda$. 
By standard dynamic programming arguments, the value function $V_k(t)$ of a representative firm satisfies the backward Hamilton--Jacobi--Bellman (HJB) equation:
\begin{equation}
    \frac{dV_k}{dt} - rV_k + \sup_{p\geq0} \left\{ \lambda(p,\bar{p},m) \left( p-\Delta V_k \right) \right\} = 0, \quad k \in \{1,\dots,K\},
    \label{eq:app_hjb}
\end{equation}
with the marginal value of inventory defined as $\Delta V_k(t) := V_k(t)-V_{k-1}(t)$. The boundary conditions are $V_k(T)=0$ and the absorbing state $V_0(t)=0$.

\paragraph{Market Intensity and Optimal Pricing.}
The baseline model utilizes a linear intensity specification:
\begin{equation}
    \lambda(p,\bar{p},m) = \left( a(m) - p + c(m)\bar{p} \right)^+,
\end{equation}
where $a(m) = 1/(1+\epsilon\eta(t))$ and $c(m) = \epsilon\eta(t)/(1+\epsilon\eta(t))$, with $\eta(t)$ representing the proportion of active firms. 
The baseline paper uses the following explicit pricing selection:
\begin{equation}
    p_k^*(t) = \frac{1}{2} \left( a(m) + c(m)\bar{p}(t) + \Delta V_k(t) \right)^+.
\end{equation}
As in the extended model, truncation of demand at zero can make the quoted
price non-unique in zero-demand states. Section~\ref{sec:extended-hjb} states
the maximizer set and the tie-breaking convention used in this paper.

\paragraph{The Kolmogorov Forward Equation.}\mbox{}\par
The aggregate market distribution $m_k(t)$ evolves deterministically according to the forward Kolmogorov equation:
\begin{equation}
    \frac{dm_k}{dt} = \lambda_{k+1}^*(t)m_{k+1}(t) - \lambda_k^*(t)m_k(t), \quad k=1,\dots,K-1,
\end{equation}
with boundary dynamics $dm_K/dt = -\lambda_K^*(t)m_K(t)$ and $dm_0/dt = \lambda_1^*(t)m_1(t)$.

\paragraph{Numerical Implementation.}
The baseline HJB--Kolmogorov system is solved numerically using explicit
time-stepping. The HJB equation is integrated backward in time using an
explicit Euler scheme, while the Kolmogorov equation is integrated forward in
time using explicit Euler. A Picard fixed-point iteration is used to resolve
the coupling between the individual optimal policy and the aggregate mean-field
variables $(\bar{p},m)$. The iteration is terminated when the supremum norm of
the difference between consecutive mean-quoted-price trajectories falls below
the tolerance
\[
\epsilon_{\mathrm{tol}}=10^{-4}.
\]
This baseline implementation serves as the numerical starting point for the
extended model analyzed in the main text.

\section{Existence of a Selected Equilibrium}
\label{app:existence}

We prove Proposition~\ref{prop:selected-equilibrium-existence}. Let
$\mathcal P(\mathcal X_{\mathrm{ext}})$ be the probability simplex. Compactness
of the action set bounds the sales rates. Together with the assumed cancellation
bound, this gives a constant $\Lambda<\infty$ bounding all exit rates. Choose
$L$ to bound the right-hand sides of all admissible forward equations and set
\begin{equation}
\begin{aligned}
\mathcal K=\bigl\{m\in C([0,T];\mathcal P(\mathcal X_{\mathrm{ext}})):
&\ m(0)=m^0,\quad \operatorname{Lip}(m)\leq L,\\
&\ \eta_m(t)\geq\eta_m(0)e^{-\Lambda t}
\text{ for all }t\in[0,T]\bigr\},
\end{aligned}
\label{eq:compact_path_set}
\end{equation}
where $\eta_m(t)=\sum_{k=-M+1}^{K}m_k(t)$. By definition, $\mathcal K$ is
nonempty, convex, closed, uniformly bounded, and equicontinuous; hence it is
compact by the Ascoli--Arzel\`a theorem.

Fix $m\in\mathcal K$. For $q\in[0,P_{\max}]$ and a vector $V$, define
\begin{equation}
 F_{V,m,t}(q)=\frac{1}{\eta_m(t)}
 \sum_{k=-M+1}^{K}m_k(t)
 \min\left\{P_{\max},
 \frac12
 \bigl(a(m(t))+c(m(t))q+\Delta V_k\bigr)^+\right\}.
 \label{eq:aggregate_price_map}
\end{equation}
The map takes $[0,P_{\max}]$ into itself. Since the positive-part and minimum
maps are one-Lipschitz,
\begin{equation}
 |F_{V,m,t}(q)-F_{V,m,t}(q')|
 \leq\frac{c(m(t))}{2}|q-q'|.
 \label{eq:aggregate_price_contraction}
\end{equation}
Since $\epsilon\geq0$, the contraction modulus is at most $1/2$, so the fixed
point $q(V,m,t)$ is unique. The lower bound on $\eta_m$, the Lipschitz continuity
of $a,c$, and the contraction estimate make $q$ continuous in $(t,V,m)$ and
locally Lipschitz in $V$, uniformly on $\mathcal K$.

For fixed $m$, insert $\bar p(t)=q(V,m,t)$ in every HJB equation and denote the
resulting right-hand side for $\dot V_k$ by $G_k(t,V;m)$. This gives the closed
ODE $\dot V=G(t,V;m)$. The preceding regularity makes $G$ continuous in time
and locally Lipschitz in $V$, while bounded prices and jump rates give
$\lVert G(t,V;m)\rVert\leq C_0+C_1\lVert V\rVert$, uniformly in $m\in\mathcal K$.
Standard ODE theory and Gronwall's inequality yield a unique, uniformly bounded
solution $V^m$ on $[0,T]$; Equation~\eqref{eq:selected_price} then gives a
continuous selected feedback $\alpha^m$ and bounded rates.

Let $\Phi(m)$ be the unique simplex-valued solution of the forward Kolmogorov
equation generated by $\alpha^m$; then $\operatorname{Lip}(\Phi(m))\leq L$.
Only a sale from $-M+1$ to $-M$ can leave the sales-active set, while a
cancellation from $-M$ re-enters it; hence
\[
 \frac{d}{dt}\eta_{\Phi(m)}(t)
 \geq-\Lambda\eta_{\Phi(m)}(t),
\]
and Gronwall's inequality yields
$\eta_{\Phi(m)}(t)\geq\eta_m(0)e^{-\Lambda t}>0$. Thus
$\Phi(\mathcal K)\subseteq\mathcal K$.

On the common bounded set containing all $V^m$, the closed HJB vector field is
uniformly Lipschitz in $V$ and continuous in $m$. ODE continuous dependence and
Equation~\eqref{eq:aggregate_price_contraction} imply that uniform convergence
$m^n\to m$ gives $V^{m^n}\to V^m$, then $\alpha^{m^n}\to\alpha^m$, and finally
$\Phi(m^n)\to\Phi(m)$. Thus $\Phi:\mathcal K\to\mathcal K$ is continuous.
Schauder's fixed-point theorem \citep{deimling1985nonlinear} therefore gives
$m^*\in\mathcal K$ such that $\Phi(m^*)=m^*$. The associated HJB solution,
selected feedback, and unique mean-quoted-price path satisfy all three
conditions in Definition~\ref{def:selected-mfg-equilibrium}. This proves the
proposition.

\bibliographystyle{plainnat}
\setlength{\bibsep}{4pt}
\bibliography{references}

\end{document}